# A new definition of the fractional Laplacian

W. Chen

Simula Research Laboratory, P. O. Box. 134, NO-1325 Lysaker, Norway

(9 September 2002)

**Keywords**: fractional Laplacian, fractional derivative, Riesz potential, Green second identity, hyper-singularity, boundary conditions

## 1. Introduction

The fractional Laplacian and the fractional derivative are two different mathematical concepts (Samko et al, 1987). Both are defined through a singular convolution integral, but the former is guaranteed to be the positive definition via the Riesz potential as the standard Laplace operator, while the latter via the Riemann-Liouville integral is not. It is noted that the fractional Laplacian can not be interpreted by the fractional derivative in the sense of either Riemann-Liouville or Caputo. Both the fractional Laplacian and the fractional derivative have found applications in many complicated engineering problems. In particular, the fractional Laplacian attracts new attentions in recent years owing to its unique capability describing anomalous diffusion problems (Hanyga, 2001).

It is, however, noted that the standard definition of the fractional Laplacian leads to a hyper-singular convolution integral and is also obscure about how to implement the boundary conditions. This purpose of this note is to introduce a new definition of the fractional Laplacian to overcome these major drawbacks. This study is carried out with the ongoing project of "mathematical and numerical modelings of medical ultrasound wave propagation" sponsored by the Simula Research Laboratory in Norway.

## 2. Riesz potential and fractional Laplacians

The fractional Laplacian is commonly considered the inverse of the Riesz potential (Gorenflo and Mainardi, 1998). The Riesz potential $I_d^s$ of order $s$ of $n$ dimensions reads (ZÄHLE, 1997; Samko et al, 1987)

$$I_d^s \varphi(x) = \frac{\Gamma[(d-s)/2]}{\pi^{s/2} 2^s \Gamma(s/2)} \int_\Omega \frac{\varphi(\xi)}{\|x-\xi\|^{d-s}} d\Omega(\xi), \quad 0<s<2, \tag{1}$$

where $\Gamma$ denotes the Euler's gamma function, and $\Omega$ is $n$-dimension integral domain. The fractional Laplacian can be defined by (e.g., see Gorenflo and Mainardi, 1998).

$$(-\Delta)_*^{s/2} \varphi(x) = -I_d^{-s} \varphi(x). \tag{2}$$

Thus, the above fractional Laplacian is also often called the Riesz fractional derivative. The above definition (2) of the fractional Laplacian can be actually restated as

$$(-\Delta)_*^{s/2} \varphi(x) = -\Delta\left[I_d^{2-s} \varphi(x)\right]. \tag{3}$$

It is known that the radial Laplacian operator has the expression

$$\Delta \varphi(x) = \frac{d^2 \varphi}{dr^2} + \frac{d-1}{r} \frac{d\varphi}{dr}, \tag{4}$$

which $r = \|x - \xi\|$. (3) can then be reduced to

$$\begin{aligned}(-\Delta)_*^{s/2} \varphi(x) &= -\frac{\Gamma[(d-2+s)/2]}{\pi^{2-s/2} 2^{2-s} \Gamma(2-s/2)} \Delta \int_\Omega \frac{\varphi(\xi)}{\|x-\xi\|^{d-2+s}} d\Omega(\xi) \\ &= -\frac{(d-2+s)s\Gamma[(d-2+s)/2]}{\pi^{2-s/2} 2^{2-s} \Gamma(2-s/2)} \int_\Omega \frac{\varphi(\xi)}{\|x-\xi\|^{d+s}} d\Omega(\xi).\end{aligned} \tag{5}$$

It is noted that (5) encounters the detrimental issues such as the hyper-singularity. An alternative way is thus presented below to define the fractional Laplacian without the perplexing issues in the Riesz fractional derivative (3)

$$(-\Delta)^{s/2}\varphi(x) = -I_d^{2-s}[\Delta\varphi(x)]$$
$$= -\frac{\Gamma[(d-2+s)/2]}{\pi^{2-s/2}2^{2-s}\Gamma(2-s/2)}\int_\Omega \frac{\Delta\varphi(\xi)}{\|x-\xi\|^{d-2+s}}d\Omega(\xi). \quad (6)$$

The Green second identity is useful to simplify (6) and can be stated as

$$\int_\Omega v\Delta\varphi d\xi = \int_\Omega \varphi\Delta v d\Omega(\xi) - \int_S \left(u\frac{\partial v}{\partial n} - v\frac{\partial \varphi}{\partial n}\right)dS(\xi), \quad (7)$$

where $S$ represents the surface of the domain, and $n$ is the unit outward normal. Let

$$v = 1/\|x-\xi\|^{d-2+s}, \quad (8)$$

and boundary conditions

$$\varphi(x)\big|_{x\in\Gamma_D} = R(x), \quad (9)$$

$$\frac{\partial\varphi(x)}{\partial n}\bigg|_{x\in\Gamma_N} = Q(x), \quad (10)$$

where $\Gamma_D$ and $\Gamma_N$ are the surface part corresponding to the Dirichlet boundary and the Neumann boundary, and using the Green second identity, the definition (6) is then reduced to

$$(-\Delta)^{s/2}\varphi(x) = -\frac{1}{h}\left\{\int_\Omega \frac{\varphi(\xi)}{\|x-\xi\|^{d+s}}d\Omega(\xi) - \right.$$

$$\left.\int_S \left[\varphi(\xi)\frac{\partial}{\partial n}\left(\frac{1}{\|x-\xi\|^{d+s}}\right) - \frac{1}{\|x-\xi\|^{d+s}}\frac{\partial \varphi(\xi)}{\partial n}\right]dS(\xi)\right\} \quad (11)$$

$$= (-\Delta)_*^{s/2}\varphi(x) + \frac{1}{h}\int_S\left[D(\xi)\frac{\partial}{\partial n}\left(\frac{1}{\|x-\xi\|^{d+s}}\right) - \frac{N(\xi)}{\|x-\xi\|^{d+s}}\right]dS(\xi),$$

where

$$h = \frac{\pi^{2-s/2} 2^{2-s} \Gamma(2-s/2)}{(d-2+s)s\Gamma[(d-2+s)/2]}. \quad (12)$$

It is seen from (11) that the presented fractional Laplacian definition is thus considered the Riesz fractional derivative (the standard fractional Laplacian) augmented with the boundary integral, which is a parallel to the fractional derivatives in the Caputo sense relative to that in the Riemann-Lioville sense. Our definition also has inherent the regularization of the hyper-singularity.

The above two definitions $(-\Delta)_*^{s/2}$ and $(-\Delta)^{s/2}$ involve in only symmetric fractional Laplacian. To clearly illustrate the basic idea of this study without loss of generality, we only consider the isotropic media in this paper. For the traditional definition of the anisotropic fractional Laplacian see Feller (1971) and Hanyga (2001). By analogy with the new definition (6) and (11), it will be straightforward to have the corresponding new expression of the anisotropic fractional Laplacian.

Albeit a long history, the research on the space fractional Laplacian still appears poor in the literature (Gorenflo and Mainardi, 1998). In recent years, some interests arise from anomalous diffusion problems. The readers are advised to find more detailed description of the fractional Laplacian from Samko et al (1987), Zaslavsky (1994), Gorenflo and Mainardi (1998), Hanyga (2001) and references therein.

## 3. FEM discretization formulation

Let the FEM discretization of a Laplacian operator be expressed as

$$-\nabla^2 p \Rightarrow K\bar{p}, \tag{13}$$

where $\bar{p}$ represents the pressure value vector at the discrete nodes, and $K$ is the positive definite FEM discretization matrix. The corresponding FEM formulation of the $s/2$ order fractional Laplacian is then obtained by

$$\left(-\nabla^2 p\right)^{s/2} \Rightarrow K^{s/2}\bar{p}, \tag{14}$$

By using a superposition analysis (Bathe and Wilson, 1976), Chen (2002) derived a FEM formulation of the power law attenuation similar to (14) in form. This reminds us that the FEM modal analysis approach demands no extra effort to solve the present fractional Laplacian model (14) which only involves the common modal parameters such as eigenvalues and eigenvectors of matrix $K$.

**References**:

Bathe, K. and Wilson, E. L (1976), *Numerical Methods in Finite Element Analysis*, Prenticle-Hall, New Jersey.

Gorenflo, R. and Mainardi, F. (1998), "Random walk models for space-fractional diffusion processes," *Fractional Calculus & Applied Analysis*, **1**, 1677-191.

Hanyga, A. (2001), "Multi-dimensional solutions of time-fractional diffusion-wave equations," to appear in *Proc. R. Soc. London A*.

Samko, S. G., Kilbas, A. A., Marichev, O. I. (1987), *Fractional Integrals and Derivatives: Theory and Applications* (Gordon and Breach Science Publishers).

ZÄHLE, M. (1997), "Fractional differentiation in the self-affine case. V - The local degree of differentiability," *Math. Nachr.* **185**, 279-306.

Zaslavsky, G. M. (1994), "Fractional kinetic equation for Hamiltonian chaos," Physica D, **76**, 110-122.